\documentclass[10pt,a4paper]{article}

\usepackage{amsmath} \usepackage{amssymb} \usepackage{amsthm}
\usepackage[dvips]{graphicx} \usepackage[all]{xy} \usepackage{eucal}
\usepackage{mathrsfs}

\SelectTips {cm}{} \CompileMatrices

\newtheorem{thm}{Theorem}[section] \newtheorem*{thm3.2*}{Theorem 3.2}

\newtheorem{lem}[thm]{Lemma} \newtheorem{prop}[thm]{Proposition}
\newtheorem*{prop*}{Proposition 2.3} 
\newtheorem{defi}[thm]{Definition} \newtheorem{rmk}[thm]{Remark}

\newenvironment{pro}[1][Proof]{{\it{#1:}} }{\hfill $\square$}
\newenvironment{pro*}[1][Proof]{{\it{#1:}} }{}

\newcommand{\mengensymb}[1]{\mathbb{#1}}
\newcommand{\A}{\mengensymb{A}} 
\newcommand{\F}{\mengensymb{F}} \newcommand{\N}{\mengensymb{N}}
\newcommand{\PP}{\mengensymb{P}} \newcommand{\Q}{\mengensymb{Q}}
 \newcommand{\Z}{\mengensymb{Z}}

\newcommand{\CC}{\mathscr{C}} \newcommand{\OO}{\mathcal{O}}
\newcommand{\DD}{\mathscr{D}} \newcommand{\XX}{\mathscr{X}}

\DeclareMathOperator{\Isom}{Isom} \DeclareMathOperator{\Out}{Out}
\DeclareMathOperator{\Aut}{Aut} \DeclareMathOperator{\Hom}{Hom}
\DeclareMathOperator{\uIsom}{\underline{\Isom}}

\newcommand{\surj}{\twoheadrightarrow}
\newcommand{\inj}{\hookrightarrow}

\newcommand{\pit}{\pi_1^{t}} \newcommand{\pitg}{\overline{\pi}_1^{t}}
\DeclareMathOperator{\EXT}{EXT} \newcommand{\pge}{{\mathscr{G}}}
\DeclareMathOperator{\Rev}{Rev}

\DeclareMathOperator{\Spec}{Spec} \DeclareMathOperator{\Sch}{Sch}
 \newcommand{\kq}{{\bar{k}}}
 \DeclareMathOperator{\FC}{FC}
\newcommand{\FCk}{\FC_k} \DeclareMathOperator{\M}{\mathscr{M}}
\newcommand{\mgn}{\M_g[n]}

\begin{document}

\title{Affine anabelian curves in positive characteristic}
\author{Jakob Stix\thanks{University of Bonn, stix@math.uni-bonn.de}}
\maketitle

\begin{abstract}
 \noindent An investigation of morphisms that coincide topologically is used to
 generalize to all characteristics and partly reprove Tamagawa's
 theorem on the Grothendieck conjecture in anabelian geometry for
 affine hyperbolic curves.  The theorem now deals with $\pi_1^{tame}$
 of curves over a finitely generated field and its effect on the sets
 of isomorphisms.
 Universal homeomorphisms are formally inverted.
\end{abstract}


\section{Introduction}

Anabelian geometry deals with anabelian categories of schemes for
which the \'etale fundamental group functor encodes many, if not all,
algebraic properties. In other words, the functor $\pi_1$ should be an
equivalence of a geometric category with a group theoretic category.

In 1983 Grothendieck announced a list of anabelian conjectures within
a letter to Faltings \cite{G83}. His ``Yoga der anabelschen
Geometrie'' declares suitable varieties - including hyperbolic curves
- over fields of absolutely finite type to be anabelian.  For a survey
and recent results see \cite{F98} and \cite{Mb}.

This paper describes a category $\FCk$ of curves over a field $k$ and
a category of pro-finite exterior Galois representations $\pge(G_k)$
which is the natural target for the tame fundamental group functor
$\pit$ applied to schemes with $k$ structure. The following theorem of
anabelian nature is proved in Section \ref{result}. It generalizes
Tamagawa's result from \cite{T97} which treats finite fields $k$ and
those of characteristic $0$.
\begin{thm3.2*} 
  Let $k$ be a field of absolutely finite type. For affine hyperbolic
  curves $C$ and $C'$ over $k$ the map
  $$\pit:\Isom_{\FCk}(C,C') \to
  \Isom_{G_k}\big(\pitg(C),\pitg(C')\big) $$
  is bijective unless both
  curves are isotrivial. In that case the map is dense injective.
\end{thm3.2*}

The method of proof relies on the idea that if a {\it fine} moduli
space $\M$ is available then isomorphisms of curves correspond to
coincidence of representing maps $\xi$.

By specialization and the finite field case $\pit$ controls the
topological component of $\xi$ for an extension of the curve over a
base $S$ of finite type over $\Z$ with function field $k$. Indeed, the
finite field case not being rigid, the $\xi$'s of two curves with
isomorphic $\pit$ differ by a ``Frobenius twist''. A priori, this
twist may vary among the various closed points of $S$. This motivates
the search for the following rigidifying result:
\begin{prop*} 
  Let $S,\M$ be of finite type over $\F_p$, $S$ be irreducible and
  reduced, and $\xi_1,\xi_2 : S \to \M$ be maps such that $\xi_1^{top}
  = \xi_2^{top}$.
  Then $\xi_1$ and $\xi_2$ differ only by a power of the Frobenius
  map.
\end{prop*}

The case of isotrivial curves needs to be dealt with seperately. They
behave differently as their $\pit$ possesses through a suitable
geometric Frobenius an automorphism of infinite order, hence $\hat{\Z}
\subset \Aut_{G_k}(\pitg)$.


\section{Topological Coincidence of Maps} \label{tcm}

For any scheme $X$ and any prime number $p$ consider $X(\F_p)$ as a
subset of the topological space underlying $X$.

\begin{lem} \label{Artin}
  Let $X$ be irreducible, of finite type over $\Spec(\Z)$, such that
  $X_\Q$ is nonempty.  Then $\bigcup_p X(\F_p)$ is Zariski-dense in
  $X$.
\end{lem}
\begin{pro}
  The generic point of $X$ is in the generic fiber $X_\Q$ which is the
  closure of its closed points. Any of these points has a number field
  as residue field and defines a closed subscheme $Z$ of dimension
  $1$.
  By the \v{C}ebotarev 
  theorem $Z$ is the closure of $\bigcup_p Z(\F_p) \subseteq \bigcup_p
  X(\F_p)$.
\end{pro}
\begin{lem}
  Let $X/k$ be irreducible and of finite type. Then any two closed
  points of $X$ lie on an irreducible curve on $X$.
\end{lem}
\begin{pro} This is a standard application of Bertini.
\end{pro}

\noindent Call the topological component of a map $f$ of schemes $f^{top}$. In
positive characteristic there is the Frobenius map $F$ which is
raising to $p^{th}$ power and has $F^{top} = id$.
\begin{prop} \label{Hauptsatz}
  Let $X,Y$ be of finite type over $\Spec(\Z)$, $X$ be irreducible and
  reduced, and $f,g : X \to Y $ be maps such that $f^{top} = g^{top}$.
  
  Then $f=g$ or $X/\F_p$. If $X/\F_p$ then $f$ and $g$ differ only by
  a power of the Frobenius map. Uniqueness of the exponent is
  equivalent to $f^{top}=g^{top}$ not being constant.
\end{prop}
\begin{pro}
  The assertion on uniqueness is clear as $f=f \circ F^m$ implies that
  the residue field at the image of the generic point of $X$ is fixed
  by $F^m$. From now on we assume that $f^{top}$ is nonconstant.
  
  In view of uniqueness we may assume $X,Y$ affine. Now the locus of
  coincidence $\{f \equiv g\}$ is closed in $X$ and contains
  $\bigcup_p X(\F_p)$. If $X_\Q \not= \emptyset$ we are done by Lemma
  \ref{Artin} .
  
  From now on assume $X/\F_p$. For $m \in \Z$ define $X_m = \{ f \circ
  F^m \equiv g\}$ or $\{f \equiv g \circ F^{-m}\}$ depending on
  the sign of $m$. Topological coincidence implies that
  $$
  X(\F_q) \subseteq \bigcup_{m \in \Z} X_m(\F_q) $$
  since
  $\F_q$-points are topologically identical if and only if they are
  $G(\F_q/\F_p)$ conjugate, and this Galois group is generated by
  Frobenius. If $q=p^r$ then it is sufficient to allow $m$ to vary
  over representatives of $\Z/r\Z$. To keep things small we choose
  representatives with minimal absolute value and thereby conserve
  symmetry:
  $$
  X(\F_{p^r}) \subseteq \bigcup_{-r/2 < m \leq r/2} X_m(\F_{p^r}) \ 
  . $$
  
  \noindent The case of $\dim X = 1$. First consider the following lemma:
\begin{lem} If $\dim X =1$ there is a constant $c$ such that for all $m$ with $X \not= X_m$ the bound $\#X_m(\F_q) \leq c \cdot p^{|m|}$ holds.
\end{lem}
\begin{pro} Choose closed immersions $X \subseteq \A^d, Y \subseteq \A^n$ and consider
  the graph $\Gamma$ of $(f,g)$
  $$
  \xymatrix@1{ X & {\Gamma} \ar[l]_{pr_1}^{\sim} } \subseteq X
  \times Y \times Y \subseteq \A^{d+2n} \subseteq \PP^{d+2n} \ . $$
  Let $(x_1, \dots,x_n,y_1, \dots,y_n)$ be coordinates of the factor
  $\A^{2n}$. For $m \in \Z$ let $R_m$ be the vanishing locus in
  $\PP^{d+2n}$ of those $n$ sections of $\mathcal{O}(p^{|m|})$
  described by $x_i^{p^m} - y_i$ or $y_i^{p^{-m}}- x_i$ (depending on
  the sign of $m$) on the affine part $\A^{d+2n}$. This $R_m$ is the
  closure of the product of $\A^d$ with the graph of $F^m$, and
  $pr_1(\Gamma \cap R_m) = X_m$.
  
  If $X_m \not= X$ then a hypersurface $H_m$ of degree $p^{|m|}$
  defined by a single suitably chosen such section suffices to cut
  down $\Gamma$ to dimension $0$.
  
  An easy intersection theoretic estimate in $\PP^{d+2n}$ with the
  closure $\bar{\Gamma}$ of $\Gamma$ gives
  $$
  \#X_m(\F_q) \leq \deg(\bar{\Gamma} \cap H_m) =
  \deg(\bar{\Gamma})\cdot \deg(H_m) = \deg(\bar{\Gamma}) \cdot p^{|m|}
  \ .$$
\end{pro}

\noindent If $X_m = X$ for no $m \in \Z$ then by the lemma
$$
\#X(\F_{p^r}) \leq \sum_{-r/2 < m \leq r/2} \! \#X_m(\F_{p^r}) \leq
c \cdot r \cdot p^{r/2} $$
for some constant $c$. This contradicts the
``Weil conjectures'' (vary $\F_{p^r}$ within the finite fields
containing the field of constants of the smooth part of $X$, then
$\#X(\F_{p^r})$ has order of magnitude $p^r$). In fact, we need only
the case of algebraic curves which is known by Hasse-Weil.

The case of $\dim X > 1$. Take $C_1,C_2 \subseteq X$ irreducible,
horizontal curves, i.e., such that $f|_{C_i}$ is not constant. By the
one dimensional case there are unique $m_i \in \Z$ such that $C_i
\subseteq X_{m_i}$. For points $x_i \in C_i$ choose an irreducible
curve $C \subset X$ passing through $x_1,x_2$. If $f(x_1) \not=
f(x_2)$ then $C\subseteq X_m$ for a unique $m\in \Z$.

If a closed point $y$ lies in $X_r \cap X_{r'}$ then $\deg(f(y)) | r -
r'$. Consequently $\deg(f(x_i))| m_i -m$. Thus
$$
\gcd\Big(\deg(f(x_1)),\deg(f(x_2))\Big) | m_1 - m_2 \ . $$
Varying
$x_1,x_2$ shows that $m_1-m_2$ has arbitrary large divisors, hence
$m_1 = m_2$. Finally $X$ is the union of such horizontal curves and so
$X = X_m$ for some $m$, as desired.
\end{pro}


\section{Preliminaries on $\pit$, Statement of Theorems} \label{pre}

\begin{defi}
  Let $S$ be a scheme. A {\bf smooth curve} over $S$ is a proper
  smooth map $p: X \to S$ of finite presentation with geometrically
  irreducible fibers of dimension $1$ together with a relative
  effective \'etale divisor $D$.
\end{defi}

Denote by $C$ the complement $X-D$ with its induced scheme structure
over $S$. Consider $C$ itself or $(X,D)$ as a short notation for the
smooth curve. The geometric genus $g$ of the fiber is a locally
constant function on $S$. The smooth curve is called affine iff
$\deg(D) > 0$ and hyperbolic iff its Euler characteristic $\chi_C = 2
- 2g - \deg(D)$ is negative.

{\bf Tame fundamental group.} Let $X$ be a normal scheme and $D
\subset X$ a divisor with normal crossing.  The pair (X,D) is
associated a tame fundamental group $\pit(X,D)$ such that the Galois
category $\Rev_{X,D}^{tame}$ of normal covers which are at most tamely
ramified along $D$ is equivalent to the category of finite continuous
$\pit(X,D)$-sets, cf.  \cite{SGA1}. Neglecting basepoints results in a
functor with values in $\pge$, the category of pro-finite groups with
exterior continuous morphisms, i.e., equivalence classes of maps up to
composition with inner automorphisms.

{\bf Exterior Galois representation.} Let $k$ always denote a field,
$\kq$ its algebraic closure. Basechange $X \times_k \kq$ is
abbreviated by $X_\kq$, and $G_k = \Aut(\kq/k)$ is the absolute Galois
group of $k$. For a smooth curve $C = (X,D)$ over $k$ the group
$\pitg(C) := \pit(X_\kq,D_\kq)$ carries a right action of $G_k$ in
$\pge$ by functorial transport of the right $G_k$-action on the
scheme. Using inverses we transform it to a left action which we
denote $\rho_C^t: G_k \to \Aut_\pge(\pitg(C))$. Clearly we obtain a
functor $\pit$ from smooth $k$-curves with values in $\pge(G_k)$, the
category of $G_k$-representations in $\pge$, i.e., pairs $(V,\rho)$
where $V \in \pge$ and $\rho:G_k \to \Aut_\pge(V) =: \Out(V)$.

{\bf ${\mathbf G_k}$ extensions.} It is known \cite{SGA1} that
the tame fundamental group of a smooth curve forms naturally an
extension
$$
1 \to \pitg(C) \to \pit(C) \to G_k \to 1 \ . $$
Maps of $k$-curves
produce morphisms between these short exact sequences with identity on
the $G_k$-part. The basepoints being neglected only
$\pitg(C)$-conjugacy classes of maps are well defined. We obtain a
functor with values in $\EXT[G_k]$, which denotes this category of
$G_k$-extensions and classes of maps. Acting on the kernel by
conjugation passes to $G_k$ if considered as an action in $\pge$, thus
a functorial construction $R : \EXT[G_k] \to \pge(G_k)$. We recover
$\rho_C^t$ from the above extension. Moreover, $R$ is an equivalence
when restricted to isomorphisms, extensions with centerfree kernel,
and $G_k$-representations on centerfree groups. The extension is
recovered by pullback: $G_k \times_{\Out(V)} \Aut(V)$ .

{\bf Fact.} If the curve is hyperbolic then $\pitg(C)$ is centerfree,
cf. \cite{F98}. The same is valid for all open subgroups as these are
$\pitg$ of covers which are hyperbolic themselves.

{\bf Topological invariance.} It is known \cite[VIII 1.1]{SGA4} that
$\pi_1$ applied to universal homeomorphisms yields isomorphisms. An
easy descente argument for tameness ensures the same behaviour for
$\pit$ of curves. The tame fundamental group is therefore not affected
by pure inseparable covers. The natural conclusion suggests to
formally invert the class of universal homeomorphisms, a task already
foreseen by
Grothendieck in his ``esquisse''. 

{\bf Frobenius} \cite[XV \S1]{SGA5}.  Fix a prime number $p$. The
Frobenius map $F$ commutes with all maps between schemes of
characteristic $p$, that is in $\Sch_{\F_p}$. If $S \in \Sch_{\F_p}$
then the diagram
$$
\xymatrix{
  X \ar@/_/[ddrr] \ar@/^/[drrr]^{F} \ar@{.>}[drr]|-{\Phi_S} \\
  & &  X(1) \ar[d] \ar[r] \ar@{}[dr]|-*+{\cdot \!\!\lrcorner}& X \ar[d] \\
  & & S \ar[r]^{F} & S } $$
defines a functor ``Frobenius twist''
$\cdot(1): Sch_S \to Sch_S$ and a natural transformation $\Phi_S :
id_{Sch_S} \to \cdot(1) $, the ``geometric $S$-Frobenius''.  They
behave well under base change: $\big(X \times_S T\big)(1) = X(1)
\times_S T$ and $\Phi_T = \Phi_S \times_S T$. The $m^{th}$ iterated
twist will be denoted by $X(m)$.

In general $X$ and its twist $X(1)$ are not isomorphic, e.g., the
twist of $\PP_k^1-\{0,1,\lambda,\infty\}$ is still genus $0$ but
punctured in $0,1,\lambda^p,\infty$.

{$\mathbf \FCk$.} Let $k$ be a field of positive characteristic.
Consider the category of smooth $k$-curves and dominant maps. Its
localization at universal homeomorphisms is easily seen to be equivalent to
its localization at geometric $k$-Frobenius maps between curves. By
Dedekind-Weber equivalence this localization can be constructed by
considering the perfection, i.e., the pure inseparable closure, of the
function field together with the unique prolongation of the set of
infinite places and maps respecting these places. Denote the resulting
category $\FCk$.

By topological invariance the tame fundamental group functor
factorizes as
$$
\pit : \FCk \to \pge(G_k), \quad C \mapsto \rho_C^t \ . $$
{\bf
  Results.} A $k$-curve $C$ is called isotrivial if $C_\kq$ is defined
over a finite field.  A field $k$ is said to be absolutely of finite
type if it is finitely generated over its prime field.  The following
theorems generalize the main result of A.~Tamagawa from \cite{T97}.
His result treats curves over base fields which are either of
characteristic $0$ or finite. Proofs will be given in the last
section.
\begin{thm} \label{resultat}
  Let $k$ be absolutely of finite type, $C$ and $C'$ be affine
  hyper\-bolic curves over $k$, such that at least one of them is not
  isotrivial. Then
  $$
  \pit:\Isom_{\FCk}(C,C') \overset{\sim}{\longrightarrow}
  \Isom_{G_k}\big(\pitg(C),\pitg(C')\big) $$
  is a bijection.
\end{thm}
\vskip -1ex
\begin{thm} \label{isotrivial} Let $k,C,C'$ be as above but $C,C'$ both isotrivial. Then the map
  $$
  \pit:\Isom_{\FCk}(C,C') \inj
  \Isom_{G_k}\big(\pitg(C),\pitg(C')\big) $$
  is injective with dense image.
\end{thm}
In particular affine hyperbolic curves $C,C'$ have 
isomorphic tame fundamental groups (with $k$-structure) if
and only if there are $m,m' \in \N$ such that $C(m) \cong C'(m')$
(with $k$-structure).

The statement of Theorem \ref{isotrivial} will become clear as one
finds an action of $\Z$, respectively $\hat{\Z}$, on the Isom-sets
compatible with the natural inclusion such that the induced map on the
quotients is bijective.  This holds essentially due to:
\begin{thm}[Tamagawa, {[Ta, 0.5]}] \label{finfield}
  Let $C,C'$ be hyperbolic curves over finite fields. Then the
  following map is a natural bijection:
  $$
  \pit : \Isom_{\Sch}\big(C,C') \overset{\sim}{\longrightarrow}
  \Isom_\pge\big(\pit(C),\pit(C')\big) \ . $$
\end{thm}
\begin{rmk}
  The condition ``affine'' in the theorem could be dropped if a
  characterization of projective hyperbolic curves over finite fields
  by their $\pit$ were available.
  
  The method of proof relies on specialization and the finite field
  case like Tamagawa's proof does, but gives a unified treatment for
  arbitrary characteristic thus also reproving the previously known.
\end{rmk}


\section{Prerequisites for the Proof}

\begin{lem}[sheaf] Let $C,C'$ be affine hyperbolic $k$-curves, and $\rho,\rho' \in \pge(G_k)$, i.e., $\rho : G_k \to \Out(V), \rho' : G_k \to \Out(V')$. Then 
\begin{itemize}
\item[(1)] $\uIsom_{\FCk}(C,C'): l/k \mapsto \Isom_{\FC_l}(C,C')$
\item[(2)] $\uIsom_{G_k}(\rho,\rho'): l/k \mapsto \Isom_{G_l}(V,V')$
\end{itemize}
are \'etale sheaves of sets on $\Spec(k)_{\acute{e}t}$. Moreover
$$
\pit: \uIsom_{\FCk}(C,C') \to
\uIsom_{G_k}\big(\rho_C^t,\rho_{C'}^t\big) $$
is a morphism of \'etale
sheaves that behaves natural w.r.t. composition.
\end{lem}
\begin{pro}
  (1) Galois descente, having localized does not matter. (2) obvious, in
  the Galois case $\uIsom_{G_k}(\rho,\rho')(k)$ are the invariants of
  the $G(l/k)$ action on $\uIsom_{G_k}(\rho,\rho')(l)$ by conjugation.
  
  The last statement is again obvious because $\pit$ is a functor and
  both Galois actions have geometric origin by conjugation with
  isomorphisms as schemes.
\end{pro}

{\bf \'Etale {$\mathbf G$}-torsors.} Let $G$ be a finite group, and
$X$ a scheme with geometric point $x$. Then almost by definition
$\Hom(\pi_1(X,x),G)$ is the set of pointed $G$-torsors $(E,e)$ on
$(X,x)$ up to isomorphism. Shifting the pointing $e \mapsto g.e$
within the fiber corresponds to composition with the inner
automorphism $g(\cdot)g^{-1}$. Hence
$$
\Hom_\pge(\pi_1(X),G) = \{ E \to X, G\text{-torsor}\} {/\cong} \ .
$$
Surjectivity is equivalent to connectedness of the torsor.

Geometrically connected tame $G$-torsors on a curve $C/k$ are
described by $\psi:\pit(C) \to G$, such that $\bar{\psi} =
\psi|_{\pitg(C)}$ is surjective. An easy diagram chase shows that
$\bar{H} = \ker(\bar{\psi})$ carries a commuting outer action of $G$
and $G_k$.
  
{\bf Definition.} Let $V \in \pge(G_k)$ be centerfree. A centerfree,
geometrically connected $G$-torsor on $V$ is a centerfree $W \in
\pge(G \times G_k)$ together with an isomorphism $G \times_{\Out(W)}
\Aut(W) \cong V$ in $\pge(G_k)$, where the $G_k$ action on the left
group is by functoriality of the construction with respect to $G$
compatible isomorphisms.

\begin{lem}
  Let $C/k$ be a hyperbolic curve.  Then we have a bijection of
  isomorphism classes:
  $$
  \left\{\begin{array}{l}
      \text{tame, geom. connected} \\
      G \text{-torsors on} \ C/k
      \end{array}
    \right\}_{/\cong} \stackrel{1:1}{\longleftrightarrow} \left\{
      \begin{array}{l}
        \text{geom. connected, center-} \\
        \text{free} \ G \text{-torsors on} \ \rho_C^t
           \end{array} 
         \right\}_{/\cong} \ . $$
\end{lem}
\begin{pro}
  $E \to C$ is mapped to $\pitg(E)$, $W \in \pge(G \times G_k)$ is
  mapped to $pr_1: (G \times G_k) \times_{\Out(W)} \Aut(W) \to G$.
  These are mutually inverse.
\end{pro}

\begin{lem}[descent] Let $C,C'$ be hyperbolic $k$-curves, $F' \to C'$ a $G$-torsor, $E' \to C'$ a geometrically connected tame $G$-torsor and $\bar{H}'$ the corresponding $G$-torsor for $\rho_{C'}^t$. Let $W'$ be a centerfree geometrically connected $G$-torsor on $V' \in \pge(G_k)$ and $V \in \pge(G_k)$. Then we have natural bijections (1), (2) and a commutative diagram (3):
\begin{itemize}
\item[(1)] $$\Big(\coprod_{G\text{-torsors}\ F/C}
  \Hom_{G,k}(F,F')\Big)_{/\cong} \stackrel{\sim}{\longrightarrow}
  \Hom_k(C,C')$$
\item[(2)] $$\Big(\coprod_{G\text{-torsors}\ W/V} \Isom_{G \times
    G_k}(W,W')\Big)_{/\cong} \stackrel{\sim}{\longrightarrow}
  \Isom_{G_k}(V,V')
  $$
\item[(3)] $$
  \xymatrix@R=0.65cm@M+0.5ex{ {\displaystyle
      \Big(\coprod_{E/C} \Isom_{G,\FCk}(E,E')\Big)_{/\cong}}
    \ar[d]_(0.6){=} \ar[r]^(0.48){\pit} &
    {\displaystyle  \Big(\coprod_{\bar{H}/\rho_C^t} \Isom_{G \times G_k}(\bar{H},\bar{H}')\Big)_{/\cong} } \ar[d]^(0.6){=} \\
    {\Isom_{\FCk}(C,C')} \ar[r]^(0.45){\pit} &
    {\Isom_{G_k}\big(\pitg(C),\pitg(C')\big) \ . }  }  $$
\end{itemize}
Here $E$ (resp. $\bar{H}$) ranges over $G$-torsors on $C$ (resp.
$\rho_C^t$), and $/\cong$ means up to equivalence induced by
isomorphisms of the ``variable'' $G$-torsors.
\end{lem}
\begin{pro} (1) $G$-torsors are $G$-quotient maps and allow pullback construction. (2) There is a map, as the ``base'' $V$ is recovered canonically by $G \times_{\Out(W)} \Aut(W)$. It is surjective by structure transport and obviously injective. (3) trivial.
\end{pro}

{\bf Specialization.}  We restate here for the convenience of the
reader Tamagawa's theorem of reconstruction of the specialization map
which is formulated entirely in grouptheoretical terms.  Base changes
are abbreviated as $X_T = X \times_S T$.

\begin{thm}[reconstruction of $sp$] \label{recspez}
  Let $K \supset R \surj k$ be a henselian discrete valuation ring and
  let $\CC/R, \CC'/R$ be smooth hyperbolic curves. Then the kernel of
  the specialization map $ sp : \pit(\CC_K) \surj \pit(\CC_k) $
  consists of the intersection of those open $H \subset \pit(\CC_K)$
  which satisfy
\begin{itemize}
\item[(i)] the image of $H$ in $G_K$ contains the inertia group $I$ of
  $R$.
\item[(ii)] the image of $I$ in $\Out(\bar{H})$ is trivial.  Here
  $\bar{H} = H \cap \pitg(\CC_K)$.
\end{itemize}
Moreover, there is a natural map $Sp$
$$
\Isom_{G_K}\big(\pitg(\CC_K),\pitg(\CC'_K)\big) \to
\Isom_{G_k}\big(\pitg(\CC_k),\pitg(\CC'_k)\big) \ . $$
\end{thm}
\begin{pro} 
  \cite[5.7]{T97}, uses known criteria for good reduction of a proper
  curve $X$ via its jacobian, minimal semistable models of $(X,D)$ and
  the combinatorics of the dual graphs for a $\Z/l\Z$-cover ramified
  along all of $D$.
\end{pro}

{\bf Level structure.} For a pro-finite group $P$ let $P^{ab}/n$
denote its maximal abelian quotient with exponent $n \in \N$. Let
$X/k$ be a proper smooth curve over $k$ of genus $g$ and $1/n \in
k^\ast$. Then $\pi_1(X_\kq)^{ab}/n$ is a $G_k$ module \'etale locally
isomorphic to $(\Z/n\Z)^{2g}$ with trivial action. A choice of an
isomorphism $\phi: (\Z/n\Z)^{2g} \cong \pi_1(X_\kq)^{ab}/n$ as $G_k$
modules is the same as equipping the curve with a level $n$ structure.
If $g\geq 2, n \geq 3$ there exists a fine moduli scheme $\mgn$
representing such pairs $(X,\phi)$, cf.  \cite[5.8,5.14]{DM}.
\begin{lem} \label{levelstructure} Let $k$ be absolutely of finite type and $C=(X,D)$ be a smooth hyperbolic curve over $k$. Then $\rho_C^t$ encodes $\pi_1(X_\kq)^{ab}/n$ as a $G_k$ module, i.e., for $C,C'$ there is a canonical map
  $$
  \Isom_{G_k}\big(\pitg(C),\pitg(C')\big) \to
  \Isom_{G_k}\big(\pi_1(X_\kq)^{ab}/n,\pi_1(X'_\kq)^{ab}/n\big) \ . $$
\end{lem}
\begin{pro} 
  There is an exact sequence of $G_k$ modules:
  $$
  0 \to \Z/n(1) \to \Z/n(1)[D(\kq)] \to \pitg(C)^{ab}/n \to
  \pi_1(X_\kq)^{ab}/n \to 0 $$
  Specialization at places of $k$
  modifies this sequence by restriction of its Galois action.  We use
  the above theorem until $k$ is finite. Now Frobenius weights
  distinguish $\pi_1(X_\kq)^{ab}/n$ as a quotient of
  $\pitg(C)^{ab}/n$.
\end{pro}

{\bf Serre rigidity.} As in the sheaf-lemma Isom's of the $G_k$
modules $\overline{\pi}_1^{ab}/n$ form a sheaf.
\begin{prop} \label{SerreRig}
  Let $X/k,X'/k$ be nonisotrivial, proper, and smooth curves of genus
  $\geq 2$, and $n \geq 3$ invertible in $k$. Then the canonical map
  $$
  \uIsom_{\FCk}(X,X') \inj
  \uIsom_{G_k}\big(\pi_1(X_\kq)^{ab}/n,\pi_1(X'_\kq)^{ab}/n\big) $$
  of
  \'etale sheaves of sets induced by $\pi_1$ is injective.
\end{prop}
\begin{pro}
  We work \'etale locally on $k$, endow $X$ with a level $n$ structure and
  obtain the characteristic map $\xi_X \in \mgn(k)$. If $f$ is a
  preimage of the identity then the corresponding $X(m) \cong X(m')$ for
  some $m,m' \in \Z$ is an isomorphism of curves with level $n$
  structure, hence $\xi_X \circ F^m = \xi_X \circ F^{m'}$. As $X$ is
  not isotrivial $m = m'$.
  
  The induced automorphism $f^\ast$ on the jacobian has finite order
  as $X(m)$ is hyperbolic and must be identity as it acts trivial on
  $n$-torsion points \cite{Serre}. In other words $f(P) - P \sim f(Q)
  - Q $ for all $P,Q \in X(\kq)$.
  
  The Lefschetz number is $\Lambda(f) = 2 -
  tr(f^\ast|H^1_{\acute{e}t}) = 2- 2g < 0$. Thus $f$ has a
  fixed point. But then $f(P) - P \sim 0$ for all $P$ and $f = id$ as
  $X$ is not $\PP^1$.
\end{pro}


\section{The Proof} \label{result}

We are going to prove now Theorem \ref{resultat}.

\begin{pro}
  Let $C=(X,D),C'=(X',D')$ be affine hyperbolic $k$-curves and $\alpha
  \in \Isom_{G_k}\big(\pitg(C),\pitg(C')\big)$. By the descente lemma
  we may prove the theorem for suitable tame covers, hence assume
  genus $g \geq 2$. (Being isotrivial holds or fails simultaneously
  for the curve and its cover). By the sheaf lemma we may enlarge $k$
  sufficiently such that the curves in question have potentially level
  $n$ structures for some $n \geq 3$.
  
  {\bf Construction of the inverse.} Choose a level structure on $X$
  and transport it to $X'$ via $\alpha$ and Lemma
  \ref{levelstructure}. Extend the data to some base $S$ of finite
  type over $\Z$ with function field $k$ and apply the following:
\begin{prop} \label{abb}
  Let $S$ be irreducible, reduced, and of finite type over $\Z$ with
  function field $k$. Consider affine hyperbolic curves $\CC =
  (\XX,\DD), \CC' = (\XX',\DD')$ over $S$ of genus $g \geq 2$ equipped
  with a level $n$ structure and generic fiber $C/k, C'/k$.
  
  If $\alpha : \pitg(C) \cong \pitg(C')$ as exterior $G_k$ modules
  s.t.  $\alpha$ respects level $n$ structures then the characteristic
  maps $\xi_\XX,\xi_{\XX'} : S \to \mgn$ representing $\XX,\XX'$
  coincide topologically.
\end{prop}
\begin{pro}
  We do induction on $\dim(S)$. We may assume $S$ normal. For
  topological coincidence it suffices to control closed points. The
  Theorem \ref{recspez} of reconstruction of specialization applied to
  the pullback of $\CC,\CC'$ over the henselization $\Spec
  \OO^h_{S,s}$ for all $s$ of codimension $1$ does the induction step.
  
  If $\dim(S) =0$ the base is a finite field $k$. We apply Theorem
  \ref{finfield} and produce an isomorphism of schemes $f: C \cong C'$
  with $\pit(f) = \alpha$ and ensure $k$ compatibility of $f$ by
  performing a suitable Frobenius-twist on say $C$. Correspondingly,
  the characteristic map is composed with a power of Frobenius and
  does not change topologically. But now $f$ is an isomorphism of
  $k$-curves with levelstructures as $\alpha=\pit(f)$ respects them,
  hence the characteristic maps coincide.
\end{pro}

Now we know that $\xi_\XX$ and $\xi_{\XX'}$ coincide topologically.
This is the point where the rigidifying effect of the result on
topological coincidence of maps develops its strength. From
Proposition \ref{Hauptsatz} we know that $\xi_\XX$ and $\xi_{\XX'}$
differ only by a unique power of Frobenius, hence there is $m \in \Z$
such that $\XX(m) \cong \XX'$ or $\XX \cong \XX'(-m)$ as $S$-curves
with level structure.

This produces an isomorphism $\tilde{\lambda}(\alpha) \in
\Isom_{\FCk}(X,X')$ which respects the effect of $\alpha$ on level $n$
structures and a natural map $\tilde{\lambda}$ in a commutative
diagram:
$$
\xymatrix@R=0.65cm@M+0.5ex{ {\Isom_{\FCk}(C,C')}
  \ar@{}[d]|{\displaystyle\bigcap} \ar[r]^(0.42){\pit} &
  {\Isom_{G_k}\big(\pitg(C),\pitg(C')\big)} \ar@{.>}@/^3ex/[l]_(0.55){\lambda}\ar@{.>}[dl]^(0.55){\tilde{\lambda}} \ar[d] \\
  {\Isom_{\FCk}(X,X')} \ar@{^(->}[r] &
  {\Isom_{G_k}\big(\pi_1(X_\kq)^{ab}/n,\pi_1(X'_\kq)^{ab}/n\big)} }
\raisebox{-10ex}{\ . } $$
Indeed, naturality is a consequence of
injectivity in the bottom line, i.e., Serre rigidity from Proposition
\ref{SerreRig} .

To detect that $\tilde{\lambda}$ even factors as $\lambda$ we apply
the descente lemma to a geometrically connected (eventually enlarge
$k$ again) tame $G$-torsor which is ramified along the whole of $D$.
The construction of $\tilde{\lambda}$ is compatible and produces a map
of $G$-torsors therefore respecting the support of ramification.

Obviously $\lambda$ is a leftinverse to $\pit$.  To see that $\pit
\circ \lambda = id$ we observe that the family of
$\lambda(\alpha|_{\alpha^{-1}(H')})$ where $H'$ varies over the open
subgroups of $\pit(C')$ defines a natural transformation $\alpha \to
\lambda(\alpha)^\ast$, compatible with $k$ structures, of functors
$\alpha,\lambda(\alpha)^\ast : \Rev^{tame}_{C'} \to \Rev^{tame}_C$
where $\lambda(\alpha)^\ast$ is pull\-back. But natural isomorphic
functors are identical on $\pit$ in $\pge(G_k)$.
\end{pro}

{\bf Isotrivial curves.} Only a sketch of proof for the case of
isotrivial curves will be given.

{\it Step 1.} If $k = \F_q, q = p^f$ then we need to quotient out the
faithfull compatible action on both Isom-sets of
$\langle F^f\rangle  \subset G_{\F_q}$
to reduce to Theorem
\ref{finfield}.

{\it Step 2.} If $C \cong C_0 \times_{\F_q} k$ such that $\F_q \subset
k$ is relatively algebraically closed, then basechange $\times_{\F_q}
k$ is an isomorphism on rational points of the Isom-scheme, which is
finite unramified, and the $G_k$ action on $\pit(C) \cong \pit(C_0)$
factors through $G_{\F_q}$. This reduces to step 1.

{\it Step 3.} We use Galois descente with some care for the density
assertion. Essentially Galois action and Frobenius commute because of
the notion of inseparable degree for morphisms in $\FCk$.


\vskip 0.3cm
{\bf Acknowledgements}
  This article contains a condensed version of the author's
  ``Di\-plom\-arbeit'', Bonn (2000), carried out under the guidance of
  his advisor Prof. Florian Pop.
  
  He deserves my sincere gratitude for invaluable support,
  encouragement, and for introducing me to the anabelian world.
  Moreover, I thank Prof. A.~Tamagawa for useful discussions at the
  MSRI and all the others who showed interest in my work. Grateful I
  enjoyed the MSRI's hospitality and excellent working conditions
  during the Galois Program in 1999 allowing substantial work for this
  article to be done.
  
  I am indebted to the referee for various remarks and comments.
  Similarly, to Eike Lau and Sabine Lange who read earlier versions of
  the manuscript I owe valuable suggestions that helped improve the
  article.


\end{document}